\documentclass[10pt]{article}

\usepackage{amsmath}
\usepackage{amsfonts}
\usepackage{amssymb}

\newtheorem{theorem}{Theorem}[section]
\newtheorem{lemma}[theorem]{Lemma}

\newtheorem{remark}[theorem]{Remark}

\newtheorem{definition}{Definition}

\def\phi{{\varphi}}

\DeclareSymbolFont{AMSb}{U}{msb}{m}{n}
\DeclareMathSymbol{\N}{\mathbin}{AMSb}{"4E}
\DeclareMathSymbol{\Z}{\mathbin}{AMSb}{"5A}
\DeclareMathSymbol{\R}{\mathbin}{AMSb}{"52}
\DeclareMathSymbol{\Q}{\mathbin}{AMSb}{"51}
\DeclareMathSymbol{\I}{\mathbin}{AMSb}{"49}
\DeclareMathSymbol{\C}{\mathbin}{AMSb}{"43}

\begin{document}
\title{The singular extremal solutions of the bilaplacian with exponential nonlinearity }
\author{  Amir Moradifam \thanks{This work is partially supported by a Killam Predoctoral
Fellowship, and is part of the author's PhD dissertation in
preparation under the supervision of N. Ghoussoub.}\\
\small Department of Mathematics,
\small University of British Columbia, \\
\small Vancouver BC Canada V6T 1Z2 \\
\small {\tt a.moradi@math.ubc.ca}
\\
%\today\\
%\date{January 20, 2005}\\
} \maketitle

\begin{abstract}
Consider the problem
\begin{eqnarray*}
\left\{ \begin{array}{ll}
\Delta^2 u=   \lambda e^{u} &\text{in } B\\
u=\frac{\partial u}{\partial n}=0 &\text{on }\partial B,
\end{array} \right.
\end{eqnarray*}
where $B$ is the unit ball in $\R^N$ and $\lambda$ is a parameter.
Unlike the Gelfand problem the natural candidate $u=-4\ln(|x|)$, for
the extremal solution, does not satisfy the boundary conditions and
hence showing the singular nature of the extremal solution in large
dimensions close to the critical dimension is challenging. D\'avila
et al. in \cite{DDGM} used a computer assisted proof to show that
the extremal solution is singular in dimensions $13\leq N\leq 31$.
Here by an improved Hardy-Rellich inequality which follows from the
recent result of Ghoussoub-Moradifam \cite{GM} we overcome this
difficulty and give a simple mathematical proof to show the extremal
solution is singular in dimensions $N\geq13$.

\end{abstract}

\section{Introduction} Consider the fourth order elliptic problem
\begin{eqnarray}\label{main.equ}
\left\{ \begin{array}{ll}
\Delta^2 u=   \lambda e^{u} &\text{in } B\\
u=\frac{\partial u}{\partial n}=0 &\text{on }\partial B,
\end{array} \right.
\end{eqnarray}
where $B$ is the unit ball in $\R^N$, $N\geq1$, $n$ is the exterior
unit normal vector and $\lambda\geq 0$ is a parameter. This problem
is the fourth order analogue of the classical Gelfand problem (see
\cite{BV}, \cite{CR}, and \cite{MP}). Recently, many authors are
intrested in fourth order equations and interesting results can be
found in \cite{AGGM}, \cite{BV}, \cite{CEGM}, \cite{DDGM},
\cite{GW}, \cite{M}, \cite{W} and the references cited therein. In
\cite {DDGM}  D\'avila et al. studied the problem (\ref{main.equ})
and showed that for each dimension $N\geq1$ there exists a
$\lambda^{*}>0$ such that for every $ 0 < \lambda < \lambda^*$,
there exists a smooth minimal (smallest) solution of
(\ref{main.equ}), while for $ \lambda
> \lambda^*$ there is no solution even in a weak sense. Moreover,
the branch $ \lambda \mapsto u_\lambda(x)$ is increasing for each $
x \in B$,  and therefore the function $u^*(x):= \lim_{\lambda
\nearrow \lambda^*} u_\lambda(x)$ can be considered as a generalized
solution that corresponds to $\lambda^*$. Now the important question
is whether $u^{*}$ is regular ($u^{*}\in L^{\infty}(B)$) or singular
($u^{*}\notin L^{\infty}(B)$). Even though there are similarities
between \ref{main.equ} and the Gelfand problem, several tools which
have been developed for the Gelfand problem, are no longer available
for (\ref{main.equ}). In \cite{DDGM} the authors developed a new
method to prove the regularity of the extremal solutions in low
dimensions and showed that for $N\leq 12$, $u^{*}$ is regular. But
unlike the Gelfand problem the natural candidate $u=-4\ln(|x|)$, for
the extremal solution, does not satisfy the boundary conditions and
hence showing the singular nature of the extremal solution in large
dimensions close to the critical dimension is challenging. D\'avila
et al. \cite{DDGM} used a computer assisted proof to show that the
extremal solution is singular in dimensions $ 13\leq N\leq31$ while
they gave a mathematical proof in dimensions $N\geq 32$. In this
paper we introduce a unified mathematical approach to deal with this
problem and show that for $N\geq13$, the extremal solution is
singular. One of our main tools is an improved Hardy-Rellich
inequality that follows from the recent result of
Ghoussoub-Moradifam about improved Hardy and Hardy-Rellich
inequalities developed in \cite{GM1} and \cite{GM}.

\section{An improved Hardy-Rellich inequality}

In this section we shall prove an improvement of classical
Hardy-Rellich inequality which will be used to prove our main result
in Section 3. It relies on the results of Ghoussoub-Moradifam in
\cite{GM} which provide necessary and sufficient conditions for such
inequalities to hold. At the heart of this characterization is the
following notion of a Bessel pair of functions.
\begin{definition} Assume that $B$ is a ball of radius $R$ in $\R^N$, $V,W \in C^{1}(0,1)$, and
$\int^{R}_{0}\frac{1}{r^{N-1}V(r)}dr=+\infty$. Say that the couple
$(V, W)$  is a {\it Bessel pair on $(0, R)$} if the ordinary
differential equation
\begin{equation*}
\hbox{ $({\rm B}_{V,W})$  \quad \quad \quad \quad \quad \quad \quad
\quad \quad \quad \quad
$y''(r)+(\frac{N-1}{r}+\frac{V_r(r)}{V(r)})y'(r)+\frac{W(r)}{V(r)}y(r)=0$
\quad \quad \quad \quad \quad \quad \quad \quad \quad \quad \quad}
\end{equation*}
has a positive solution on the interval $(0, R)$.
\end{definition}

\begin{theorem}\label{GHR} ({\bf Ghoussoub-Moradifam} \cite{GM}) Let $V$ and $W$ be positive radial $C^1$-functions   on $B\backslash \{0\}$, %%@
where $B$ is a ball centered at zero with radius $R$ in $\R^N$ ($N \geq 1$) such that  %%@
$\int^{R}_{0}\frac{1}{r^{N-1}V(r)}dr=+\infty$ and $\int^{R}_{0}r^{N-1}V(r)dr<+\infty$. The %%@
following statements are then equivalent:

\begin{enumerate}

\item $(V, W)$ is a Bessel pair on $(0, R)$ and $\beta (V, W; R) \geq 1$.

\item $ \int_{B}V(x)|\nabla u |^{2}dx \geq \int_{B} W(x)u^2dx$ for all $u \in C_{0}^{\infty}(B)$.

\item If $\lim_{r \rightarrow 0}r^{\alpha}V(r)=0$ for some $\alpha<
N-2$ and $W(r)-\frac{2V(r)}{r^2}+\frac{2V_r(r)}{r}-V_{rr}(r)\geq 0$ on $(0, R)$, %%@
then the above are equivalent to
\[
\hbox{$\int_{B}V(x)|\Delta u |^{2}dx \geq  \int_{B} W(x)|\nabla  %%@
u|^{2}dx+(N-1)\int_{B}(\frac{V(x)}{|x|^2}-\frac{V_r(|x|)}{|x|})|\nabla u|^2dx$ \quad  for all  $u %%@
\in C^{\infty}_{0}(B)$.}
\]
\end{enumerate}

\end{theorem}

As an application we have the following improvement of the classical
Hardy-Rellich inequality.

\begin{theorem}
Let $N\geq 5$ and $B$ be the unit ball in $\R^N$. Then the following
improved Hardy-Rellich inequality holds for all $u \in
C^{\infty}_{0}(B)$.
\begin{equation}\label{HR1}
\int_{B}|\Delta u|^2\geq
\frac{(N-2)^2(N-4)^2}{16}\int_{B}\frac{u^2}{(|x|^2-\frac{9}{10}|x|^{\frac{N}{2}+1})(|x|^2-|x|^{\frac{N}{2}})}+\frac{(N-1)(N-4)^2}{4}
\int_{B}\frac{u^2}{|x|^2(|x|^2-|x|^{\frac{N}{2}})}.
\end{equation}
As a consequence the following improvement of classical
Hardy-Rellich inequality holds:
\begin{equation}\label{HR2}
\int_{B}|\Delta u|^2\geq
\frac{N^2(N-4)^2}{16}\int_{B}\frac{u^2}{|x|^2(|x|^2-|x|^{\frac{N}{2}})}.
\end{equation}
\end{theorem}
{\bf Proof.} Let $\varphi:=r^{-\frac{N}{2}+1}-\frac{9}{10}$. Since
\begin{equation*}
-\frac{\varphi''+\frac{(N-1)}{r}\varphi'}{\varphi}=\frac{(N-2)^2}{4}.\frac{1}{r^2-\frac{9}{10}r^{\frac{N}{2}+1}},
\end{equation*}
$(1,\frac{(N-2)^2}{4}\frac{1}{r^2-\frac{9}{10}r^{\frac{N}{2}+1}})$
is a bessel pair on $(0,1)$. By Theorem \ref{GHR} the following
inequality holds for all $u \in C_{0}^{\infty}(B)$.
\begin{equation}\label{AHR}
\int_{B}|\Delta u|^2 dx\geq \frac{(N-2)^2}{4} \int_{B}\frac{|\nabla
u|^2}{|x|^2-\frac{9}{10}|x|^{\frac{N}{2}+1}}+(N-1)\int_{B}\frac{|\nabla
u|^2}{|x|^2}.
\end{equation}
Let $V(r):=\frac{1}{r^2-\frac{9}{10}r^{\frac{N}{2}+1}}$. Then
\begin{equation}
\frac{V_{r}}{V}=-\frac{2}{r}+\frac{9}{10}(\frac{N-2}{2})\frac{r^{\frac{N}{2}-2}}{1-\frac{9}{10}r^{\frac{N}{2}-1}}\geq-\frac{2}{r},
\end{equation}
and $\psi(r)=r^{-\frac{N}{2}+2}-1$ is a positive super-solution for
the ODE
\begin{equation}\label{ODE}
y''+(\frac{N-1}{r}+\frac{V_{r}}{V})y'(r)+\frac{W_{1}(r)}{V(r)}y=0,
\end{equation}
where
\[W_{1}(r)=\frac{(N-4)^2}{4(r^2-r^{\frac{N}{2}})(r^2-\frac{9}{10}r^{\frac{N}{2}+1})}.\]

Hence the ODE (\ref{ODE}) has actually a positive solution and by
Theorem \ref{GHR} we have
\begin{equation}
\int_{B}\frac{|\nabla
u|^2}{|x|^2-\frac{9}{10}|x|^{\frac{N}{2}+1}}\geq(\frac{N-4}{2})^2\int_{B}\frac{u^2}{(|x|^2-\frac{9}{10}|x|^{\frac{N}{2}+1})(|x|^2-|x|^{\frac{N}{2}})}.
\end{equation}
Similarly
\begin{equation}
\int_{B}\frac{|\nabla u|^2}{|x|^2}\geq
(\frac{N-4}{2})^2\int_{B}\frac{u^2}{|x|^2(|x|^2-|x|^{\frac{N}{2}})}.
\end{equation}
Combining the above two inequalities with (\ref{AHR}) we get
(\ref{HR1}). $\Box$

\section{Main results}
In this section we shall prove that the extremal solution $u^{*}$ of
the problem (\ref{main.equ}) is singular in dimensions $N\geq13$.
The next lemma will be our main tool to guarantee that $u^{*}$ is
singular for $N\geq13$. The proof is based on an upper estimate by a
singular stable sub-solution.

\begin{lemma}\label{sing-lem}
Suppose there exist $\lambda'>0$ and a radial function $u \in
H^{2}(B)\cap W^{4,\infty}_{loc}(B\setminus \{0\})$ such that
\begin{equation}\label{cond1}
\Delta^2 u\leq \lambda'e^{u} \ \ \hbox{for all} \ \ 0<r<1,
\end{equation}

\begin{equation}
u(1)=0, \ \ \ \ \frac{\partial u}{\partial n}(1)=0,
\end{equation}
\begin{equation}
u \notin L^{\infty}(B),
\end{equation}
and
\begin{equation}\label{cond2}
\beta \int_{B}e^{u} \varphi^2\leq \int_{B}(\Delta \varphi)^2\ \
\hbox{for all} \ \ \varphi \in C^{\infty}_{0}(B),
\end{equation}

for some $\beta > \lambda'$. Then $u^{*}$ is singular and
\begin{equation}\label{est1}
\lambda^{*} \leq \lambda'
\end{equation}
\end{lemma}
{\bf Proof.} By Lemma 2.6 in \cite{DDGM} we have (\ref{est1}).
Define
\begin{equation}
\alpha:=\ln(\frac{\lambda'}{\gamma \lambda^*}),
\end{equation}
where $\frac{\lambda'}{\beta}<\gamma<1$ and let $\bar{u}:=u+\alpha$.
We claim that
\begin{equation}\label{claim}
u^{*}\leq \bar{u} \ \ \hbox{in} \ \ B.
\end{equation}
To prove this, we shall show that for $\lambda<\lambda^{*}$
\begin{equation}
u_{\lambda}\leq \bar{u} \ \ \hbox{in} \ \ B.
\end{equation}
Indeed, we have
\begin{eqnarray*}
\Delta^2(\bar{u})=\Delta^2(u)\leq
\lambda'e^{u}=\lambda'e^{-\alpha}e^{\bar{u}}=\gamma \lambda^{*}
e^{\bar{u}}.
\end{eqnarray*}
To prove (\ref{claim}) it suffices to prove it for
$\gamma\lambda^{*}<\lambda<\lambda^{*}$. Fix such $\lambda$ and
assume that $(\ref{claim})$ is not true. Let
\begin{equation*}
R_{1}:=\sup\{0\leq R\leq 1\mid u_{\lambda}(R)=\bar{u}(R)\}.
\end{equation*}
Since $\bar{u}(1)=\alpha>0=u_{\lambda}(1)$, we have $0<R_{1}<1$,
$u_{\lambda}(R_{1})=\bar{u}(R_{1})$, and $u'_{\lambda}(R_{1})\leq
\bar{u}'(R_{1})$. Now consider the following problem
\begin{eqnarray*}
\left\{ \begin{array}{ll}
\Delta^2 u=   \lambda e^{u} &\text{in } \Omega\\
u=u_{\lambda}(R_{1}) &\text{on }\partial \Omega\\
\frac{\partial u}{\partial n}=u'_{\lambda}(R_{1}) &\text{on
}\partial \Omega.
\end{array} \right.
\end{eqnarray*}
Obviously $u_{\lambda}$ is a solution to the above problem while
$\bar{u}$ is a sub-solution to the same problem. Moreover $\bar{u}$
is stable since,
\[\lambda<\lambda^{*}\]
and hence
\[\lambda e^{\bar{u}}\leq \lambda^{*}e^{\alpha}e^{u}= \frac{\lambda'}{\gamma} e^{u}<\beta e^{u}.\]
We deduce $\bar{u}\leq u_{\lambda}$ in $B_{R_{1}}$ which is
impossible, since $\bar{u}$ is singular while $u_{\lambda}$ is
smooth. This establishes (\ref{claim}). From (\ref{claim}) and the
above two inequalities we have
\[\lambda^{*}e^{u^{*}}\leq \lambda^{*}e^{a}e^{u}=\frac{\lambda'}{\gamma}e^{u}.\]
Since $\frac{\lambda'}{\gamma}<\beta$,
\[\inf_{\varphi \in C^{\infty}_{0}}(B)\frac{\int_{B}(\Delta \varphi)^2-\lambda^{*}e^{u^{*}}}{\int_{B}\varphi^2}>0.\]
This is not possible if $u^{*}$ is a smooth solution. $\Box$ \\

In the following, for each dimension $N\geq13$, we shall construct
$u$ satisfying all the assumptions of Lemma \ref{sing-lem}. Define
\[w_{m}:=-4\ln(r)-\frac{4}{m}+\frac{4}{m}r^{m}, \ \ m>0,\]
and let $H_{N}:=\frac{N^2(N-4)^2}{16}$. Now we are ready to prove
our main result.

\begin{theorem}\label{main.sing}  The following upper bounds on $\lambda^*$ hold in large dimensions.
\begin{enumerate}
\item If $N \ge 32$, then Lemma \ref{sing-lem} holds with $u:=w_2$, $\lambda_{N}'=8(N-2)(N-4)e^2$ and $\beta=H_{N}>\lambda'_{N}$.
\item If $13\le N \le 31$, then Lemma \ref{sing-lem} holds with $u:=w_{3.5}$ and $\lambda'_N<\beta_N$ given in Table \ref{table:summary}.
\end{enumerate}
The extremal solution is therefore singular for dimensions $N\geq
13$.

\end{theorem}

{\bf Proof.} 1) Assume first that $N\geq 32$, then
\[8(N-2)(N-4)e^{2}<\frac{N^2(N-4)^2}{16},\]
and
\[\Delta^2 w_{2}=8(N-2)(N-4)\frac{1}{r^4}\leq8(N-2)(N-4)e^2e^{w_{2}}.\]
Moreover,
\[8(N-2)(N-4)e^2\int_{B}e^{w_{2}}\varphi^2\leq H_{n}\int_{B}e^{-4\ln(|x|)}\varphi^2=H_{n}\int_{B}\frac{\varphi^2}{|x|^2}\leq \int_{B}|\Delta \varphi|^2.\]
Thus it follows from Lemma \ref{sing-lem} that $u^{*}$ is singular
and $\lambda^{*}\leq 8(N-2)(N-4)e^2$.

2) Assume $13 \leq N\leq 31$. We shall show that $u=w_{3.5}$
satisfies the assumptions of Lemma \ref{sing-lem} for each dimension
$13\leq N\leq 31$. Using Maple, for each dimension $13\leq N\leq
31$, one can verify that  inequality (\ref{cond1}) holds for
$\lambda'_{N}$ given by Table \ref{table:summary}. Then, by using
Maple again, we show that there exists $\beta_{N}>\lambda'_{N}$ such
that
\begin{eqnarray*}
\frac{(N-2)^2(N-4)^2}{16}\frac{1}{(|x|^2-0.9|x|^{\frac{N}{2}+1})(|x|^2-|x|^{\frac{N}{2}})}&+&\frac{(N-1)(N-4)^2}{4}
\frac{1}{|x|^2(|x|^2-|x|^{\frac{N}{2}})}\\ &\geq&
\beta_{N}e^{w_{3.5}}.
\end{eqnarray*}
The above inequality and improved Hardy-Rellich inequality
(\ref{HR1}) guarantee that the stability condition (\ref{cond2})
holds for $\beta_{N}>\lambda'$. Hence by Lemma \ref{sing-lem} the
extremal solution is singular for $13\leq N\leq 31$. The values of
$\lambda_{N}$ and $\beta_{N}$ are shown in Table
\ref{table:summary}.

\begin{table}[ht]
\caption{Summary} % title of Table
\centering % used for centering table
\begin{tabular}{c c c } % centered columns (4 columns)
\hline\hline %inserts double horizontal lines
N & $\lambda'_{N}$ & $\beta_{N}$  \\ [0.5ex] % inserts table
%heading
\hline % inserts single horizontal line
 $N\geq 32$ & $8(N-2)(N-4)e^2$ & $H_{n}$  \\
31 & 20000 & 86900  \\
30 & 18500 & 76500  \\
29 & 17000 & 67100 \\
28 & 16000 & 58500 \\
27 & 14500 & 50800 \\
26 & 13500 & 43870  \\
25 & 12200 & 37630  \\
24 & 11100 & 32050  \\
23 & 10100 & 27100  \\
22 & 9050 & 22730 \\
21 & 8150 & 18890  \\
20 & 7250 & 15540  \\
19 & 6400 & 12645  \\
18 & 5650 & 10155  \\
17 & 4900 & 8035  \\
16 & 4230 & 6250  \\
15 & 3610 & 4765  \\
14 & 3050 & 3545  \\
13 & 2525 & 2560  \\[1ex] % [1ex] adds vertical space
\hline %inserts single line
\end{tabular}
\label{table:summary} % is used to refer this table in the text
\end{table}

\begin{remark} The values of $\lambda'_{N}$ and $\beta_{N}$ in Table
\ref{table:summary} are not optimal.
\end{remark}

\begin{remark}
The improved Hardy-Rellich inequality (\ref{HR1}) is crucial to
prove that $u^{*}$ is singular in dimensions  $N\geq 13$. Indeed by
the classical Hardy-Rellich inequality and $u:=w_{3.5}$, Lemma
\ref{sing-lem} only implies that $u^{*}$ is singular in dimensions
$N\geq22$.
\end{remark}

{\bf Acknowledgment:} I would like to thank Professor Nassif
Ghoussoub, my supervisor, for his valuable suggestions, constant
support, and encouragement.

\end{document}